\documentclass{amsart}
\usepackage{amssymb}
\usepackage{amsmath}
\usepackage{amsfonts}

\newtheorem{theorem}{Theorem}
\theoremstyle{plain}

\newtheorem{exercise}{Exercise}

\numberwithin{equation}{section}
\input{tcilatex}

\begin{document}
\title[Orbital Measures]{Existence and Uniqueness of Orbital\ Measures}
\author{Michael Barnsley}
\address{Australian National University, Canberra\\
ACT, Austalia}
\email{Mbarnsley@aol.com}
\date{ \today }

\begin{abstract}
We note an elementary proof of the existence and uniqueness of a solution $%
\mu \in \mathbb{P(X)}$ to the equation $\mu =p\mu _{0}+q\widehat{F}\mu $.
Here $\mathbb{X}$ is a topological space, $\mathbb{P(X)}$ is the set of
Borel measures of unit mass on $\mathbb{X}$, $\mu _{0}\in $ $\mathbb{P(X)}$
is given, $p>0$, and $q\geq 0$ with $p+q=1$. The transformation $\widehat{F}:%
\mathbb{P(X)\rightarrow P(X)}$ is defined by $\widehat{F}\upsilon
=\tsum\limits_{n=1}^{N}p_{n}\upsilon \circ f_{n}^{-1}$ where $f_{n}:\mathbb{%
X\rightarrow X}$ is continuous, $p_{n}>0$ for $n=1,2,...,N$, $N$ is a finite
strictly positive integer, and $\tsum\limits_{n=1}^{N}p_{n}=1$. This problem
occurs in connection with iterated function systems (IFS).
\end{abstract}

\maketitle

\section{Introduction}

In this note we present an elementary proof of the existence and uniqueness
of a solution to the pair of equations 
\begin{equation}
\mu =p\mu _{0}+q\widehat{F}\mu \text{ with }\mu \in \mathbb{P(X)}\text{.}
\label{nonhomequation}
\end{equation}%
Here $\mathbb{X}$ is a topological space, $\mathbb{P(X)}$ is the set of
normalized positive Borel measures on $\mathbb{X}$, $\mu _{0}\in $ $\mathbb{%
P(X)}$ is given, $p>0$, $q\geq 0$, and $p+q=1$. The transformation $\widehat{%
F}:\mathbb{P(X)\rightarrow P(X)}$ is defined by 
\begin{equation}
\widehat{F}\upsilon =\tsum\limits_{n=1}^{N}p_{n}\upsilon \circ f_{n}^{-1}%
\text{ for all }\upsilon \in \mathbb{P(X)}\text{,}  \label{operatorequation}
\end{equation}%
where $f_{n}:\mathbb{X\rightarrow X}$ is continuous, $p_{n}>0$ for $%
n=1,2,...,N$, $N$ is a finite strictly positive integer and $%
\tsum\limits_{n=1}^{N}p_{n}=1$. Equation (\ref{nonhomequation}) occurs in
connection with iterated function systems (IFS). In \cite{BaDe} extra
restictions are placed on the functions $f_{n}$ and the space $\mathbb{X}$
to enable a contraction mapping argument to be used, similar to the one used
in \cite{Hu81} for the case $p=0$.

Given $f:\mathbb{X\rightarrow X}$ we use the same symbol $f$ to denote $f:%
\mathbb{P(X)\rightarrow P(X)}$ where $f(\upsilon )=\upsilon \circ f^{-1}$
for all $\upsilon \in \mathbb{P(X)}$.

Let $\Omega _{\{1,2,...,N\}}^{\prime }$ denote the set of all finite strings
of symbols taken from the set $\{1,2,...,N\}$. That is, $\sigma \in \Omega
_{\{1,2,...,N\}}^{\prime }$ iff there is $K\in \mathbb{N}$ and $\sigma
=\sigma _{1}\sigma _{2}...\sigma _{K}$ where $\sigma _{k}\in \{1,2,...,N\}$
for each $k\in \{1,2,...,K\}$.

\begin{theorem}
\label{measureorbittheorem} Let the setting defined above be true. Then the
measure $\mu \in \mathbb{P(X)}$ which is defined by%
\begin{equation}
\mu =p\mu _{0}+\tsum\limits_{\sigma \in \Omega _{\{1,2,...,N\}}^{\prime
},\left\vert \sigma \right\vert \geq 1}pq^{\left\vert \sigma \right\vert
}p_{\sigma _{1}}p_{\sigma _{2}}...p_{\sigma _{\left\vert \sigma \right\vert
}}f_{\sigma _{1}}\circ f_{\sigma _{2}}\circ ...\circ f_{\sigma _{\left\vert
\sigma \right\vert }}(\mu _{0})\text{.}  \label{orbitmeasureequation}
\end{equation}
is the unique solution to Equation (\ref{nonhomequation}).
\end{theorem}

\begin{proof}
Let $\mathcal{B}(\mathbb{X)}$ denote the set of Borel subsets of $\mathbb{X}$%
. The series 
\begin{equation*}
p\mu _{0}(B)+\tsum\limits_{\sigma \in \Omega _{\{1,2,...,N\}}^{\prime
},\left\vert \sigma \right\vert \geq 1}pq^{\left\vert \sigma \right\vert
}p_{\sigma _{1}}p_{\sigma _{2}}...p_{\sigma _{\left\vert \sigma \right\vert
}}f_{\sigma _{1}}\circ f_{\sigma _{2}}\circ ...\circ f_{\sigma _{\left\vert
\sigma \right\vert }}(\mu _{0})(B)
\end{equation*}%
is absolutely convergent, uniformly in $B\in \mathcal{B}(\mathbb{X)}$,
because it consists of non-negative terms and is bounded above,
term-by-term, by the absolutely convergent series%
\begin{equation*}
p+\tsum\limits_{\sigma \in \Omega _{\{1,2,...,N\}}^{\prime },\left\vert
\sigma \right\vert \geq 1}pq^{\left\vert \sigma \right\vert }p_{\sigma
_{1}}p_{\sigma _{2}}...p_{\sigma _{\left\vert \sigma \right\vert }}=1\text{.}
\end{equation*}
Hence the value $\mu (B)$ is well-defined for all $B\in \mathcal{B}(\mathbb{%
X)}$ and $\mu :\mathcal{B}(\mathbb{X)\rightarrow \lbrack }0,1]$. Notice that 
$\mu (\mathbb{X})=1$.

Let us define 
\begin{equation*}
\rho _{0}=\mu _{0}\text{ and }\rho _{n}=\tsum\limits_{\sigma \in \Omega
_{\{1,2,...,N\}}^{\prime },\left\vert \sigma \right\vert =n}p_{\sigma
_{1}}p_{\sigma _{2}}...p_{\sigma _{\left\vert \sigma \right\vert }}f_{\sigma
_{1}}\circ f_{\sigma _{2}}\circ ...\circ f_{\sigma _{\left\vert \sigma
\right\vert }}(\mu _{0})\text{ for }n=1,2,..\text{.}
\end{equation*}%
Then it is readily verified that $\rho _{n}\in \mathbb{P(X)}$ and we can
rewrite Equation (\ref{orbitmeasureequation}) as 
\begin{equation*}
\mu =\tsum\limits_{n=0}^{\infty }pq^{n}\rho _{n}\text{.}
\end{equation*}%
Let $\{\mathcal{O}_{m}\in \mathcal{B}(\mathbb{X}):m=1,2,...\}$ be a sequence
such that $\bigcup\limits_{m=1}^{\infty }\mathcal{O}_{m}\in \mathcal{B}(%
\mathbb{X})$ and 
\begin{equation*}
\mathcal{O}_{m_{1}}\cap \mathcal{O}_{m_{2}}=\varnothing 
\end{equation*}%
for all $m_{1},m_{2}\in \mathbb{N}$ with $m_{1}\neq m_{2}$. Then%
\begin{equation*}
\tsum\limits_{m=0}^{\infty }\mu (\mathcal{O}_{m})=\tsum\limits_{m=1}^{\infty
}\tsum\limits_{n=0}^{\infty }pq^{n}\rho _{n}(\mathcal{O}_{m})\text{.}
\end{equation*}%
Since the series on the right is absolutely convergent, we can interchange
the order in which the two summations are evaluated, which yields%
\begin{equation*}
\tsum\limits_{m=1}^{\infty }\mu (\mathcal{O}_{m})=\tsum\limits_{n=0}^{\infty
}pq^{n}\tsum\limits_{m=1}^{\infty }\rho _{n}(\mathcal{O}_{m})=\tsum%
\limits_{n=0}^{\infty }pq^{n}\rho _{n}(\bigcup\limits_{m=1}^{\infty }%
\mathcal{O}_{n})=\mu (\bigcup\limits_{m=1}^{\infty }\mathcal{O}_{n})\text{.}
\end{equation*}%
It follows that $\mu $ is indeed a measure on $\mathcal{B}(\mathbb{X)}$ and,
since $\mu (\mathbb{X})=1$, it follows that $\mu \in \mathbb{P(X)}$.

In order to prove the second part of the theorem we note that, since all of
the series involved are absolutely convergent, it suffices to demonstate
that the algebra works out correctly, term-by-term. Substituting from
Equation (\ref{orbitmeasureequation}) into Equation (\ref{nonhomequation})
we find%
\begin{eqnarray*}
&&\text{r.h.s. of Equation (\ref{nonhomequation})} \\
&=&p\mu _{0}+q\tsum\limits_{n=1}^{N}p_{n}f_{n}(p\mu
_{0}+\tsum\limits_{\sigma \in \Omega _{\{1,2,...,N\}}^{\prime },\left\vert
\sigma \right\vert \geq 1}pq^{\left\vert \sigma \right\vert }p_{\sigma
_{1}}p_{\sigma _{2}}...p_{\sigma _{\left\vert \sigma \right\vert }}f_{\sigma
_{1}}\circ f_{\sigma _{2}}\circ ...\circ f_{\sigma _{\left\vert \sigma
\right\vert }}(\mu _{0})) \\
&=&p\mu _{0}+q\tsum\limits_{n=1}^{N}p_{n}f_{n}(p\mu
_{0}+\tsum\limits_{m=1}^{N}qp_{0}p_{m}f_{m}(\mu _{0}) \\
&&+\tsum\limits_{\sigma \in \Omega _{\{1,2,...,N\}}^{\prime },\left\vert
\sigma \right\vert \geq 2}pq^{\left\vert \sigma \right\vert }p_{\sigma
_{1}}p_{\sigma _{2}}...p_{\sigma _{\left\vert \sigma \right\vert }}f_{\sigma
_{1}}\circ f_{\sigma _{2}}\circ ...\circ f_{\sigma _{\left\vert \sigma
\right\vert }}(\mu _{0})) \\
&=&p_{0}\mu _{0}+\tsum\limits_{n=1}^{N}pqp_{n}f_{n}(\mu
_{0})+\tsum\limits_{n=1}^{N}\tsum%
\limits_{m=1}^{N}pq^{2}p_{n}p_{m}f_{n}(f_{m}(\mu _{0}))+ \\
&&\tsum\limits_{n=1}^{N}\tsum\limits_{\sigma \in \Omega
_{\{1,2,...,N\}}^{\prime },\left\vert \sigma \right\vert \geq
2}pq^{\left\vert \sigma \right\vert +1}p_{n}p_{\sigma _{1}}p_{\sigma
_{2}}...p_{\sigma _{\left\vert \sigma \right\vert }}f_{\sigma _{1}}\circ
f_{\sigma _{2}}\circ ...\circ f_{\sigma _{\left\vert \sigma \right\vert
}}(\mu _{0})) \\
&=&p_{0}\mu _{0}+\tsum\limits_{\sigma \in \Omega _{\{1,2,...,N\}}^{\prime
},\left\vert \sigma \right\vert \geq 1}p_{0}(1-p_{0})^{\left\vert \sigma
\right\vert }p_{\sigma _{1}}p_{\sigma _{2}}...p_{\sigma _{\left\vert \sigma
\right\vert }}f_{\sigma _{1}}\circ f_{\sigma _{2}}\circ ...\circ f_{\sigma
_{\left\vert \sigma \right\vert }}(\mu _{0}) \\
&=&\text{l.h.s. of Equation (\ref{nonhomequation}).}
\end{eqnarray*}%
In order to prove uniqueness, suppose that $\upsilon \in \mathbb{P(X)}$
obeys 
\begin{equation*}
\upsilon =p\mu _{0}+q(p_{1}f_{1}(\upsilon )+p_{2}f_{2}(\upsilon
)+...+p_{N}f_{N}(\upsilon ))\text{.}
\end{equation*}%
Then, by repeatedly substituting from the left-hand-side into the
right-hand-side we find that $\upsilon $ can be represented by the same
absolutely convergent series as $\mu $, whence $\upsilon =\mu $.
\end{proof}

In \cite{barnsley06}, Chapter 3, we refer to the unique solution of Equation
(\ref{nonhomequation}) as the \textbf{orbital measure} associated with the
IFS $\{\mathbb{X};f_{1},f_{2},...,f_{N};p_{1},p_{2},...,p_{N}\}$, the
probabilities $p$ and $q$, and the \textbf{condensation measure} $\mu _{0}$.
Note that we have\textit{\ not} required that the underlying space be
complete or compact, or that the IFS be contractive or contractive on the
average.

Notice that the expressions above could have been written down and handled
more succinctly in terms of the operator $\widehat{F}:\mathbb{%
P(X)\rightarrow }\mathbb{P(X)}$ defined by Equation (\ref{operatorequation}%
). $\widehat{F}$ acts linearly on the space of linear combinations of Borel
measures on $\mathbb{X}$. Using this notation the series expansion in
Equation (\ref{orbitmeasureequation}) can be written as%
\begin{eqnarray*}
\mu  &=&p(1-q\widehat{F})^{-1}\mu _{0} \\
&=&p\tsum\limits_{m=0}^{\infty }q^{m}(\widehat{F})^{m}\mu _{0}\text{.}
\end{eqnarray*}

\begin{exercise}
Let $\mathbb{X=}[0,1)\subset \mathbb{R}$ with the usual topology. Let $%
S_{\{f\}}(\mathbb{X)}$ be the semigroup generated by the function $%
f:[0,1)\rightarrow \lbrack 0,1)$ defined by $f(x)=\frac{1}{2}+\frac{x}{2}$.
Let $\mu _{0}\in \mathbb{P([}0,1))$ denote a normalized Borel measure all of
whose mass is contained in $[0,\frac{1}{2})$. That is, $\mu _{0}([0,\frac{1}{%
2}))=1,$ and $\mu _{0}((\frac{1}{2},1))=0$. Then the associated orbital
measure $\mu \in \mathbb{P(X)}$ satisfies%
\begin{equation*}
\mu =p\mu _{0}+qf(\upsilon )=\sum_{n=0}^{\infty }pq^{n}f^{\circ n}(\mu _{0})%
\text{.}
\end{equation*}%
What happens as $p\rightarrow 0$? Do we get a solution to $\upsilon
=f(\upsilon )$ with $\upsilon \in \mathbb{P(X)}$? Show that for each $x\in
\lbrack 0,1)$ we have 
\begin{equation*}
\lim_{p_{0}\rightarrow 0}\mu ([0,x])=0\text{.}
\end{equation*}%
Conclude that we do not obtain, in the limit, a solution to $\upsilon
=f(\upsilon )$ with $\upsilon \in \mathbb{P(X)}$. What happens if the
interval $[0,1)$ is replaced by $[0,1]$?
\end{exercise}

\end{document}